# CORRECTION
# IMPROPER REGULAR CONDITIONAL DISTRIBUTIONS

By Teddy Seidenfeld, Mark J. Schervish and Joseph B. Kadane

*Carnegie Mellon University*

A strict inequality appears in Definition 6 where a weak inequality is needed. We reproduce Definition 6 here.

DEFINITION 6. Fix $\omega$ and consider those $A$ such that $\omega \in A \in \mathcal{A}$. If for some $\omega \in A \in \mathcal{A}$, $P(A|\mathcal{A})(\omega) = 0$, say that $P(\cdot|\mathcal{A})$ is *maximally improper at* $\omega$. Otherwise, if for each $\omega \in A \in \mathcal{A}$, $1 \geq P(A|\mathcal{A})(\omega) > 0$, say that the rcd is *modestly proper at* $\omega$.

At the bottom of page 1614, we are not precise in the definition of a Borel space. The condition should have read that there is a one-to-one measurable function with measurable inverse between $(\Omega, \mathcal{B})$ and $(E, \mathcal{E})$, where $E$ is a Borel subset of the reals and $\mathcal{E}$ is the Borel $\sigma$-field of subsets of $E$. After the remaining corrections below, our use of the term "Borel space" conforms with this definition.

Some conditions were left out of Theorem 4 and Lemma 3. The proof of Lemma 3 also had some errors that made it almost impossible to follow. Finally, the proof of Theorem 4 was said to be straightforward from Theorem 3. We include here the restatements of both results with the missing conditions, the revised proof of Lemma 3, and a proof of Lemma 4. The only application of Lemma 4 given in the original paper is to the proof of Corollary 2. The additional conditions given here are satisfied in that case.

THEOREM 4. *Assume that $\mathcal{A}$ is an atomic sub-$\sigma$-field of $\mathcal{B}$. Let $(\Theta, \mathcal{D})$ be a Borel space, with a probability measure $\mu$. For each $\theta \in \Theta$, let $P_\theta$ be a probability on $\mathcal{B}$ such that for every $B \in \mathcal{B}$, $P_\theta(B)$ is a $\mathcal{D}$-measurable function of $\theta$. Let $P(\cdot)$ be defined on $\mathcal{B}$ by $P(\cdot) = \int_\Theta P_\theta(\cdot) \, d\mu(\theta)$. Assume that, for $\mu$-almost all $\theta$, $P_\theta(\cdot|\mathcal{A})$ is a maximally improper rcd for $P_\theta$ and* *that it is $\mathcal{A} \otimes \mathcal{D}$-measurable as a function of $(\omega, \theta)$. Also, assume that the set*

$$B^* = \{(\omega, \theta) : P_\theta(\cdot|\mathcal{A}) \text{ is maximally improper at } \omega\},$$

*is in $\mathcal{A} \otimes \mathcal{D}$. Then there is a maximally improper version of $P(\cdot|\mathcal{A})$.*

LEMMA 3. *Let $(\Theta, \mathcal{D})$ be a Borel space, with a probability measure $\mu$. For each $\theta \in \Theta$, let $P_\theta$ be a probability on $\mathcal{B}$ such that for every $B \in \mathcal{B}$, $P_\theta(B)$ is a $\mathcal{D}$-measurable function of $\theta$. Define the probability $P$ on $\mathcal{B}$ by $P(B) = \int_\Theta P_\theta(B) \, d\mu(\theta)$. Let $\mathcal{A}$ be a sub-$\sigma$-field of $\mathcal{B}$. Also, let $P_\theta(\cdot|\mathcal{A})$ denote an rcd for each $P_\theta$ that is $\mathcal{A} \otimes \mathcal{D}$-measurable as a function of $(\omega, \theta)$. Then, for each $\omega$ there exists a probability $\nu_\omega$ on $\mathcal{D}$ such that for all $B \in \mathcal{B}$*

$$\text{(1)} \qquad \int_\Theta P_\theta(B|\mathcal{A})(\omega) \, d\nu_\omega(\theta)$$

*is a version of $P(B|\mathcal{A})$. Also, these versions form an rcd.*

PROOF. Let $\mathcal{E}$ be the product $\sigma$-field $\mathcal{B} \otimes \mathcal{D}$. For each $E \in \mathcal{E}$, define

$$E_\theta = \{\omega : (\omega, \theta) \in E\},$$
$$E^\omega = \{\theta : (\omega, \theta) \in E\},$$

the $\theta$- and $\omega$-sections of $E$. Standard arguments like those of Billingsley ([1], Section 18) allow us to conclude that $E_\theta \in \mathcal{B}$ for all $\theta$, and $P_\theta(E_\theta)$ is a $\mathcal{D}$-measurable function of $\theta$. Define

$$Q(E) = \int_\Theta P_\theta(E_\theta) \, d\mu(\theta),$$

which is easily seen to be a probability on $\mathcal{E}$. Let $\pi_1(\omega, \theta) = \omega$ and $\pi_2(\omega, \theta) = \theta$ be the coordinate projections, which are $\mathcal{E}$-measurable. Let $\mathcal{A}' = \pi_1^{-1}(\mathcal{A})$ and $\mathcal{D}' = \pi_2^{-1}(\mathcal{D})$, which are sub-$\sigma$-fields of $\mathcal{E}$. Every $\mathcal{A}'$-measurable function must be an $\mathcal{A}$-measurable function of $\pi_1$. Because $(\Theta, \mathcal{D})$ is a Borel space, there exists an rcd for $\pi_2$ given $\mathcal{A}'$ relative to $Q$, $Q(\cdot|\mathcal{A}')$. We will denote $Q(\pi_2^{-1}(D)|\mathcal{A}')(\omega, \theta)$ by $\nu_\omega(D)$. In similar fashion to the arguments earlier in the proof, $\nu_\omega(E^\omega)$ is $\mathcal{A}$-measurable as a function of $\omega$ for all $E \in \mathcal{E}$. Define

$$Q_0(E) = \int \nu_\omega(E^\omega) \, dP(\omega).$$

For each $A \in \mathcal{A}$ and $D \in \mathcal{D}$, we have

$$Q_0(A \times D) = \int I_A \nu_\omega(D) \, dP(\omega) = Q(A \times D).$$

It follows that $Q_0 = Q$ on all of $\mathcal{A} \otimes \mathcal{D}$.

For each $\omega$, (1) is a probability. We need to show that it is $\mathcal{A}$-measurable as a function of $\omega$. We have assumed that $P_\theta(\cdot|\mathcal{A})(\omega)$ is $\mathcal{A} \otimes \mathcal{D}$ measurable, so we can approximate it from below by a sequence $\{\phi_n\}_{n=1}^\infty$ of nonnegative simple functions. In similar fashion to the argument at the beginning of this proof, $\nu_\omega(E^\omega)$ is $\mathcal{A}$-measurable for all $E \in \mathcal{A} \otimes \mathcal{D}$. It follows that $\int \phi_n(\omega, \theta) \, d\nu_\omega(\theta)$ is $\mathcal{A}$-measurable for each $n$, and (1) is a limit of $\mathcal{A}$-measurable functions.

To complete the proof, we show that, for each $A \in \mathcal{A}$ and $\mathcal{B} \in \mathcal{B}$, the integral of (1) over $A$ equals $P(A \cap B)$:

$$\int_A \int_\Theta P_\theta(B|\mathcal{A})(\omega) \, d\nu_\omega(\theta) \, dP(\omega) = \int I_A(\omega) P_\theta(B|\mathcal{A})(\omega) \, dQ_0(\omega, \theta)$$

$$= \int I_A(\omega) P_\theta(B|\mathcal{A})(\omega) \, dQ(\omega, \theta)$$

$$= \int \int I_A(\omega) P_\theta(B|\mathcal{A})(\omega) \, dP_\theta(\omega) \, d\mu(\theta)$$

$$= \int P_\theta(A \cap B) \, d\mu(\theta) = P(A \cap B),$$

where the first equality is from the definition of $Q_0$, the second follows from the fact that $Q_0 = Q$ on $\mathcal{A} \otimes \mathcal{D}$, the third is from the definition of $Q$, the fourth is from the definition of $P_\theta(\cdot|\mathcal{A})$ and the last is the meaning of $P_\theta$. □

PROOF OF THEOREM 4. Because $\mathcal{A}$ is atomic, $P_\theta(\cdot|\mathcal{A})$ is maximally improper at $\omega$ if and only if $P_\theta(a(\omega)|\mathcal{A})(\omega) = 0$, where $a(\omega)$ is the $\mathcal{A}$-atom containing $\omega$. Hence, we can rewrite the set $B^*$ as

$$B^* = \{(\omega, \theta) : P_\theta(a(\omega)|\mathcal{A})(\omega) = 0\},$$

whose $\theta$-sections satisfy

$$B^*_\theta = \{\omega : P_\theta(a(\omega)|\mathcal{A})(\omega) = 0\} \in \mathcal{B}.$$

For each $\theta$ such that $P_\theta(\cdot|\mathcal{A})$ is maximally improper, $B^*_\theta$ has inner $P_\theta$ measure 1. Hence $P_\theta(B^*_\theta) = 1$, a.e. $[\mu]$. By standard arguments, $P_\theta(B^*_\theta)$ is $\mathcal{D}$-measurable, and it follows that

$$Q(B^*) = \int_\Theta P_\theta(B^*_\theta) \, d\mu(\theta) = 1,$$

where $Q$ was constructed in the proof of Lemma 3.

Similarly, the $\omega$-sections of $B^*$ satisfy

$$B^{*\omega} = \{\theta : P_\theta(a(\omega)|\mathcal{A})(\omega) = 0\} \in \mathcal{D}.$$

For each $\omega$, let $\nu_\omega$ be the measure from Lemma 3. Then $\nu_\omega(B^{*\omega})$ is $\mathcal{D}$-measurable. Since $B^* \in \mathcal{A} \otimes \mathcal{D}$, we have

$$1 = Q(B^*) = Q_0(B^*) = \int_\Omega \nu_\omega(B^{*\omega}) \, dP(\omega),$$

where $Q_0$ was constructed in the proof of Lemma 3. So, there is a set $C \in \mathcal{B}$ with $P(C) = 1$ and for all $\omega \in C$, $\nu_\omega(B^{*\omega}) = 1$. It follows that, for each $\omega \in C$, there is a set $E(\omega) \in \mathcal{D}$ with $\nu_\omega(E(\omega)) = 1$ such that $P_\theta(a(\omega)|\mathcal{A})(\omega) = 0$ for all $\theta \in E(\omega)$. Let $P(\cdot|\mathcal{A})$ be the version guaranteed by Lemma 3. Then, for each $\omega \in C$,

$$P(a(\omega)|\mathcal{A})(\omega) = \int_\Theta P_\theta(a(\omega)|\mathcal{A})(\omega) \, d\nu_\omega(\theta) = 0.$$

This means that $P(\cdot|\mathcal{A})$ is maximally improper. $\square$

**Acknowledgments.** The authors would like to thank P. Berti and P. Rigo both for bringing these points to their attention and for helping to patch the proofs. They would also like to thank a referee of the correction for an extremely careful reading that further helped correct and simplify the proofs.

T. SEIDENFELD  
DEPARTMENT OF PHILOSOPHY  
CARNEGIE MELLON UNIVERSITY  
PITTSBURGH, PENNSYLVANIA 15213  
USA  
E-MAIL: teddy@stat.cmu.edu

M. J. SCHERVISH  
J. B. KADANE  
DEPARTMENT OF STATISTICS  
CARNEGIE MELLON UNIVERSITY  
PITTSBURGH, PENNSYLVANIA 15213  
USA  
E-MAIL: mark@stat.cmu.edu  
kadane@stat.cmu.edu